arXiv:1706.08732v1 [math.OC] 27 Jun 2017# On efficiently solving the subproblems of a level-set method for fused lasso problems

Xudong Li[*], Defeng Sun[†] and Kim-Chuan Toh[‡]June 27, 2017

**Abstract**

In applying the level-set method developed in [Van den Berg and Friedlander, SIAM J. on Scientific Computing, 31 (2008), pp. 890–912 and SIAM J. on Optimization, 21 (2011), pp. 1201–1229] to solve the fused lasso problems, one needs to solve a sequence of regularized least squares subproblems. In order to make the level-set method practical, we develop a highly efficient inexact semismooth Newton based augmented Lagrangian method for solving these subproblems. The efficiency of our approach is based on several ingredients that constitute the main contributions of this paper. Firstly, an explicit formula for constructing the generalized Jacobian of the proximal mapping of the fused lasso regularizer is derived. Secondly, the special structure of the generalized Jacobian is carefully extracted and analyzed for the efficient implementation of the semismooth Newton method. Finally, numerical results, including the comparison between our approach and several state-of-the-art solvers, on real data sets are presented to demonstrate the high efficiency and robustness of our proposed algorithm in solving challenging large-scale fused lasso problems.

## 1 Introduction

Let $p: \Re^n \to \Re$ be the fused lasso regularizer, i.e,

$$p(x) := \lambda_1 \|x\|_1 + \lambda_2 \|Bx\|_1, \quad \forall x \in \Re^n,$$

where $\lambda_1, \lambda_2 \geq 0$ are given parameters, $B \in \Re^{(n-1) \times n}$ is the matrix defined by $Bx = [x_1 - x_2, x_2 - x_3, \ldots, x_{n-1} - x_n]^T$, $\forall x \in \Re^n$. First proposed in [41], the fused lasso regularizer is designed to encourage the sparsity in both the coefficients and their successive differences. This regularizer is particularly suitable for problems with features that can be ordered in some meaningful ways. In this paper, we consider the following fused lasso problem:

$$\min \left\{ p(x) \mid \|Ax - b\| \leq \delta \right\}, \tag{1}$$

where $A \in \Re^{m \times n}$ is a given matrix, $b \in \Re^m$ and $\delta > 0$ are given data. Comparing to the following regularzied least squares form

$$\min \left\{ \frac{1}{2} \|Ax - b\|^2 + p(x) \right\}, \tag{2}$$

the least squares constrained formulation (1) is widely believed to be computationally more challenging because of the complicated geometry of the feasible region. Yet, formulation (1) is sometimes preferred in real-data modeling since one can always control the noise level of the model through tuning the acceptable tolerance – the parameter $\delta$.

---

[*]Department of Mathematics, National University of Singapore, 10 Lower Kent Ridge Road, Singapore 119076 (matlixu@nus.edu.sg).

[†]Department of Mathematics and Risk Management Institute, National University of Singapore, 10 Lower Kent Ridge Road, Singapore 119076 (matsundf@nus.edu.sg).

[‡]Department of Mathematics and Institute of Operations Research and Analytics, National University of Singapore, 10 Lower Kent Ridge Road, Singapore 119076 (mattohkc@nus.edu.sg).
1

A potentially feasible approach for solving problem (1) is the recently developed level-set method [42, 43, 1]. It has been shown to possess superior performance in many interesting least squares constrained optimization problems including the basis pursuit denoising [42, 43] and matrix completion [43]; see [43, 1] for more examples. When applied to problem (1), the level-set method developed in [42, 43, 1] executes an iterative root finding procedure for solving the following univariate nonlinear equation

$$\phi(\tau) = \delta,$$

where $\phi$ is the value function of the level-set minimization problem resulting from exchanging the objective and constraint functions in problem (1), i.e.,

$$\phi(\tau) := \min \left\{ \|Ax - b\| \mid p(x) \leq \tau \right\}, \quad \tau \geq 0. \tag{3}$$

Therefore, instead of solving problem (1) directly, the level-set method solves a sequence of minimization problems of form (3) which are parameterized by $\tau$. As noted already in [42], this approach depends critically on the availability of an efficient solver for problem (3). Note that the algorithms proposed in [42, 43, 1] for problem (3) require a closed-form representation or an efficient computation of the metric projector over the feasible set $F := \{x \in \Re^n \mid p(x) \leq \tau\}$. However, due to the composite structure in $p$, no efficient approach for computing the metric projector is currently available. Fortunately, as one will see shortly, the level-set method can be carefully designed to avoid the potentially highly expensive computations of the metric projector. Of course, it is an interesting topic to develop an efficient way to compute the metric projector, but we will leave it for future research and would not focus on this issue in this paper.

Although the computation of the projector mentioned above is hindered by the composite structure in $p$, the proximal mapping of $p$ in fact can be computed in a fairly easy manner. Indeed, Friedman et al. in [14] showed that the proximal mapping of $p$ can be obtained, in a semi-closed-form expression, through the composition of the proximal mappings of two individual parts of $p$, i.e., the $\ell_1$ norm $\|\cdot\|$ and the TV-norm $\|B \cdot \|$. This decomposition property has been further studied in [44] and is termed as "prox-decomposition". From [44], one can see that many interesting regularizers, such as the elastic-net regularizer [47], Berhu regularizer [28] and many feature grouping regularizers, enjoy this special "prox-decomposition" property. Although the metric projectors on the level set of the aforementioned regularizers are difficult to calculate, the exceptional "prox-decomposition" feature can be exploited to design fast methods to compute their proximal mappings.

The "prox-decomposition" property of the fused lasso regularizer, together with the difficulties of computing the metric projector $\Pi_F$, implies that when the level-set method is applied to solve the fused lasso problem (1), one should solve a sequence of regularized least squares problems. More specifically, we show that the level-set method is based on iteratively solving the following nonlinear equation

$$\varphi(\mu) := \|Ax^* - b\| = \delta, \tag{4}$$

where $x^*$ is an optimal solution of the following regularized least squares problem

$$\min \left\{ \frac{1}{2} \|Ax - b\|^2 + \mu p(x) \right\} \tag{5}$$

and $\mu \geq 0$ is a varying parameter. Indeed, careful analyses on the properties of $\varphi$ are conducted to make the above procedure executable. Our approach here sheds new light on how the level-set method shall be used for solving least squares constrained optimization problems in the form of problem (1) when the regularizer $p$ possesses complicated yet special structures. As the backbone of the level-set method, in this paper, we aim to provide a highly efficient solver for solving the above fused lasso regularized least squares problem (5).

To achieve the goal above, we propose to use the semismooth Newton augmented Lagrangian (Ssnal) method to solve problem (5). Here we are motivated by the fact that Ssnal has already proven its superior performance in solving the $\ell_1$ regularized least squares problems [21]. Note that since the objective in (5) is convex piecewise liner-quadratic, the asymptotically superlinear convergence of Ssnal has been shown in [21, Theorem 3]. With the guarantee of this fast local convergence of the augmented Lagrangian method (ALM), the sole key challenge to obtain a fast practical algorithm is in designing a highly efficient semismooth Newton method for solving the subproblem at each ALM iteration. To this end, a computable generalized Jacobian of the proximal mapping of the fused lasso regularizer $p$ is critically needed. However, such a generalized Jacobian is not available in the literature possibly due to the presence of the TV-norm and the complicated composite structure in $p$. Fortunately, under the "prox-decomposition" property and the tools for analyzing the generalized Jacobian of the metric projector over a polyhedral set developed in [18, 22], we are able to derive a nontrivial formula for constructing the generalized Jacobian of the proximal mapping of the fused lasso regularizer. Just as in [21], we need to carefully extract and analyse the special second order sparsity structure in the generalized Jacobian to ensure the efficient implementation of the semismooth Newton method. In particular, based on the second order sparsity structure, we also design sophisticated numerical techniques to efficiently solve the large scale linear systems involved in the semismooth Newton method. As the reader may expect, our Ssnal is highly efficient and robust and it substantially outperforms several state-of-the-art solvers in solving large scale fused lasso problems with real data.

The remaining parts of this paper are organized as follows. As a preliminary, the next section is devoted to studying the properties of the value function $\varphi$ and the generalized Jacobian of the solution mapping of a strongly convex quadratic programming (QP) problem with parameters. In Section 3, the explicit formula for constructing the generalized Jacobian of the proximal mapping of the fused lasso regularizer is derived. The sparsity structure of the generalized Jacobian is also carefully extracted and analyzed. Section 4 focuses on using the Ssnal to solve the regularized least squares subproblems. Efficient numerical techniques for implementing the Ssnal are also discussed. In Section 5, we conduct extensive numerical experiments with large-scale real data to evaluate the performance of SSNAL in solving various fused lasso problems. We conclude our paper in the final section.

Before we move to the next section, here we list some notation to be used later. For any given vector $y \in \Re^n$, we denote by $\mathrm{Diag}(y)$ the diagonal matrix whose $i$th diagonal element is $y_i$. For any given matrix $A \in \Re^{m \times n}$, we use $\mathrm{Ran}(A)$ and $\mathrm{Null}(A)$ to denote the range space and null space of $A$, respectively. We use $I_n$ to denote the $n$ by $n$ identity matrix in $\Re^{n \times n}$ and $N^\dagger$ to denote the Moore-Penrose pseudo-inverse of a given matrix $N \in \Re^{m \times n}$. Similarly, $\mathbf{O}_n$ and $\mathbf{E}_n$ are used to denote the $n \times n$ zero matrix and the $n \times n$ matrix of all ones, respectively. Given any index set $\alpha \subseteq \{1, \ldots, n\}$, we denote its cardinality by $|\alpha|$. For any given proper closed convex function $q : \Re^n \to (-\infty, +\infty]$, the proximal mapping $\mathrm{Prox}_q(\cdot)$ of $q$ is defined by

$$\mathrm{Prox}_q(x) = \arg\min_{z \in \Re^n} \left\{ q(z) + \frac{1}{2}\|z - x\|^2 \right\}, \quad x \in \Re^n.$$

We will often make use of the following Moreau identity $\mathrm{Prox}_{tq}(x) + t\mathrm{Prox}_{q^*/t}(x/t) = x$, where $t > 0$ is a given parameter, and $q^*$ is the Fenchel conjugate function of $q$. See [34, Section 31] for more discussions on proximal mappings.

## 2 Some preliminary results

In this section, we shall first analyze the properties of the value function $\varphi$ to make the level-set method executable and then study the generalized Jacobian of the solution mapping of a strongly convex QP with parameters, which forms the foundation for calculating the generalized Jacobian of the fused lasso regularizer.

### 2.1 Properties of the value function $\varphi$

For the purpose to study the properties of the value function $\varphi$ defined in (4), we consider a similar problem to (5) but with a more general regularizer

$$\min \left\{ \frac{1}{2} \|Ax - b\|^2 + \mu \kappa(x) \right\}, \tag{6}$$

where $\kappa : \Re^n \to (-\infty, +\infty]$ is a nonnegative positively homogeneous convex function such that $\kappa(0) = 0$, i.e., a gauge function [34, Section 15]. Here, we further assume that $\kappa$ is a convex piecewise linear function, i.e., a polyhedral convex function ([34, Section 19] and [37, Theorem 2.49]). Obviously, the fused lasso regularizer $p$ is a special instance of $\kappa$. From [5], one can observe that piecewise linear gauge functions are extremely important in handling some ill-posed inverse problems. The polar of $\kappa$ is defined by

$$\kappa^\circ(y) := \inf\{u \geq 0 \mid \langle x, y \rangle \leq u\kappa(x), \forall x \in \Re^n\}.$$

Note that $\kappa^\circ$ is also a guage function [34, Theorem 15.1]. It is not difficult to prove that

$$\kappa^\circ = \delta^*_{\kappa \leq 1} \quad \& \quad \kappa^* = \delta^*_{\kappa^\circ \leq 1},$$

e.g., see [13, Proposition 2.1 (iii) and (iv)]. Since $\kappa$ is a polyhedral convex function, its level set $\{x \in \Re^n \mid \kappa(x) \leq 1\}$ is obviously a polyhedral convex set. Then, from [34, Corollary 19.2.1], we know that $\kappa^\circ$ is also a polyhedral convex function. Now we can write the dual of problem (6) as follows:

$$\max\left\{-\frac{1}{2}\|\xi\|^2 + \langle b, \xi\rangle \mid \kappa^\circ(A^T\xi) \leq \mu\right\}. \tag{7}$$

For every $\mu$, let $\Omega(\mu)$ and $\xi(\mu)$ be the solution set of the primal problem (6) and the dual problem (7), respectively. Obviously, as (multi-)functions of $\mu$, $\text{dom}(\Omega) = \text{dom}(\xi) = \{\mu \in \Re \mid \mu \geq 0\}$. It is also not difficult to see that $\xi$ is a single-valued mapping on its domain. After all these preparations, we have the following proposition.

**Proposition 1.** *It holds that*

(i) *for any $\mu \geq 0$, $\xi(\mu) = b - Ax$, $\forall x \in \Omega(\mu)$, i.e., $\|b - Ax\|$ is invariant on the solution set $\Omega(\mu)$;*

(ii) *for all $\mu \geq \mu_\infty := \kappa^\circ(A^Tb)$, $\xi(\mu) = b$ and $0 \in \Omega(\mu)$;*

(iii) *$\xi$ is a piecewise affine function on $\text{dom}(\xi)$;*

(iv) *if $\mu_\infty > 0$, then for any $0 \leq \mu_1 < \mu_2 \leq \mu_\infty$, $\|\xi(\mu_1)\| < \|\xi(\mu_2)\|$, i.e., $\|b - Ax_1\| < \|b - Ax_2\|$ for all $x_1 \in \Omega(\mu_1)$ and $x_2 \in \Omega(\mu_2)$.*

*Proof.* (i) The equation follows directly from the KKT condition corresponding to problems (6) and (7).

(ii) Obviously, for any $\mu \geq 0$, we have $\frac{1}{2}\|b\|^2 \geq -\frac{1}{2}\|\xi(\mu)\|^2 + \langle b, \xi(\mu)\rangle$. Hence, when $\mu \geq \mu_\infty$, $b$ is the unique optimal solution to (7). From (i), we have for $\mu \geq \mu_\infty$ that $\kappa(x_\mu) = 0$ and $0 \in \Omega(\mu)$.

(iii) Let $S(\mu) := \{\xi \in \Re^m \mid \kappa^\circ(A^T\xi) \leq \mu\}$. Obviously, $\text{graph}\,S := \{(\mu, \xi) \in \Re_+ \times \Re^m \mid \kappa^\circ(A^T\xi) \leq \mu\} = \text{epi}(\kappa^\circ A^T)$. Since $\kappa^\circ$ is polyhedral convex, we know from [34, Corollary 19.3.1] that $\kappa^\circ A^T$ is a polyhedral convex function and thus $\text{epi}(\kappa^\circ A^T)$ is a polyhedral convex set. Therefore, $S$ is a graph-convex polyhedral multifunction [33, Section 2]. Then, it follows from [33, Proposition 2.4] that $\xi$ is a polyhedral multifunction. Since $\xi$ is a single-valued mapping on its domain, from [37, Excercise 2.48] and [12, Excercise 5.6.14], we know that $\xi$ is piecewise affine on its domain.

(iv) It is easy to see that $\|\xi(\mu)\|$ is a nondecreasing function of $\mu \geq 0$, e.g., see [4, Lemma 9.2.1]. We prove the required result by contradiction. Suppose that there exist $0 \leq \mu_1 < \mu_2 \leq \mu_\infty$, such that $\|\xi(\mu_1)\| = \|\xi(\mu_2)\|$. Then, $\|\xi(\mu)\| = \|\xi(\mu_1)\|$ for all $\mu \in [\mu_1, \mu_2]$. Since $\xi$ is a piecewise affine function on $[\mu_1, \mu_2]$, we have that $\xi(\mu_1) = \xi(\mu_2)$.

Now $\kappa^\circ(A^T\xi(\mu_2)) = \kappa^\circ(A^T\xi(\mu_1)) \leq \mu_1 < \mu_2$. Thus the constraint $\kappa^\circ(A^T\xi) \leq \mu_2$ is inactive at the solution $\xi(\mu_2)$, and we easily get $\xi(\mu_2) = b$ from the optimality condition. From here we have $\mu_2 > \mu_1 \geq \kappa^\circ(A^T\xi(\mu_2)) = \kappa^\circ(A^Tb) = \mu_\infty$, which contradicts the fact that $\mu_2 \leq \mu_\infty$. This completes the proof of the proposition. □

**Remark 1.** *The piecewise affine property of $\xi$ implies that $\|\xi(\mu)\|$ is piecewise smooth as well as piecewise convex on $\mu \geq 0$. Indeed, on each piece, $\xi$ can be represented as $\xi(\mu) = \alpha\mu + \beta$ for some given vectors $\alpha, \beta \in \Re^m$ with $\|\xi(\mu)\| = \sqrt{\|\beta\|^2 + 2\langle\alpha, \beta\rangle\mu + \|\alpha\|^2\mu^2}$.*

From Proposition 1, we know that the value function $\varphi$ given in (4) is well-defined and nondecreasing. In particular, it is strictly increasing on $\mu \in [0, p^\circ(A^Tb)]$, where $p^\circ$ is the polar of the fused lasso regularizer. This monotonicity and the boundedness of $\varphi(\mu)$ naturally imply that the bisection method or the secant method can be employed to solve the univariate nonlinear equation (4), i.e., $\varphi(\mu) = \delta$, where $\delta > 0$ is the given parameter in problem (1). In fact, from Remark 1, [1, Theorem A.1] and [31, Thereom 3.2], we can prove that the secant method converges at least Q-superlinearly when certain mild nondegeneracy conditions are satisfied. Under the assumption that the inequality constraint $\|Ax - b\| \leq \delta$ is active at any optimal solution of problem (1), we know that $x_{\mu^*} \in \Omega(\mu^*)$ is an optimal solution to problem (1), where $\mu^*$ is a solution of the nonlinear equation (4).

## 2.2 The generalized Jacobian of the solution mapping of a strongly convex QP

In this section, we study the generalized Jacobian of the solution mapping of a parametric strongly convex QP. The results presented here will form the foundation for studying the generalized Jacobian of the proximal mapping of the fused lasso regularizer.

Consider a nonempty polyhderal convex set $\mathcal{D} \subseteq \Re^n$ expressed in the following form

$$\mathcal{D} := \{x \in \Re^n \mid Cx \geq c,\ Dx = d\},$$

where $C \in \Re^{k \times n}$ and $D \in \Re^{l \times n}$ are given matrices, $c \in \Re^k$ and $d \in \Re^l$ are given vectors. Without loss of generality, we assume that $\mathrm{rank}(D) = l$, $l \leq n$. Given a point $x \in \Re^n$, consider the solution mapping of the following strongly convex quadratic programming (QP) problem

$$s(x) := \mathrm{argmin}\left\{\frac{1}{2}\langle s,\ Qs\rangle - \langle x,\ s\rangle \mid s \in \mathcal{D}\right\}, \tag{8}$$

where $Q \in \Re^{n \times n}$ is a given symmetric and positive definite matrix.

Given the strong convexity of the objective in problem (8), since $\mathcal{D} \neq \emptyset$, the solution mapping $s(\cdot)$ is well-defined and single-valued throughout $\Re^n$. When $Q = I_n$, the above QP reduces to the metric projection problem and $s$ is exactly the projector over $\mathcal{D}$. Similarly, $s(x)$ can be viewed as a skewed projector of $x$ onto the polyhedral set $\mathcal{D}$ in the case when $Q \neq I_n$. By using [12, Proposition 4.14] and the change-of-variables technique, we can easily show that $s$ is piecewise affine on $\Re^n$. Based on this property, one may further use the results of Pang and Ralph [30] to characterize the B-subdifferential and the corresponding Clarke generalized Jacobian [8] of $s$. However, the calculations of these generalized Jacobians can be a very difficult task to accomplish numerically for an arbitrary polyhedral set $\mathcal{D}$ and a general positive definite matrix $Q$. To circumvent this difficulty for the case with $Q = I_n$, Han and Sun in [18] defined a computable generalized Jacobian of the metric projector over $\mathcal{D}$. More recently, Han and Sun's generalized Jacobian has been further studied and used for developing efficient algorithms for solving QP problems with Birkhoff polytope constraints [22]. Here, we aim to extend Han and Sun's computable generalized Jacobian, which is defined for the metric projector only, to the solution mapping of a strongly convex QP.

From the definition of $s(x)$, we know that there exist multipliers $\lambda \in \Re^k$ and $\mu \in \Re^l$ such that

$$\begin{cases} Qs(x) - x + C^T\lambda + D^T\mu = 0, \\ Cs(x) - c \geq 0, \quad Ds(x) - d = 0, \\ \lambda \leq 0, \quad \lambda^T(Cs(x) - c) = 0. \end{cases} \tag{9}$$

Let $M(x)$ be the set of multipliers associated with $x$, i.e.,

$$M(x) := \{(\lambda, \mu) \in \Re^m \times \Re^p \mid (x, \lambda, \mu)\ \text{satisfies (9)}\}.$$

Since $M(x)$ is a nonempty polyhedral convex set containing no lines, it has at least one extreme point denoted as $(\bar{\lambda}, \bar{\mu})$ [34, Corollary 18.5.3]. Denote
$$I(x) := \{i \mid C_i s(x) = c_i, i = 1, \ldots, m\}, \tag{10}$$
where $C_i$ is the $i$th row of the matrix $C$. Define a collection of index sets:
$$\mathcal{K}(x) := \{ K \subseteq \{1, \ldots, m\} \mid \exists (\lambda, \mu) \in M(x) \text{ s.t. } \operatorname{supp}(\lambda) \subseteq K \subseteq I(x),$$
$$[C_K^T \ D^T] \text{ is of full column rank}\},$$
where $\operatorname{supp}(\lambda)$ denotes the support of $\lambda$, i.e., the set of indices $i$ such that $\lambda_i \neq 0$ and $C_K$ is the matrix consisting of the rows of $C$, indexed by $K$. As noted in [18], the set $\mathcal{K}(x)$ is nonempty due to the existence of the extreme point $(\bar{\lambda}, \bar{\mu})$ of $M(x)$. Define the following multi-valued mapping $\mathcal{P} : \Re^n \rightrightarrows \Re^{n \times n}$:
$$\mathcal{P}(x) := \left\{ P \in \Re^{n \times n} \mid P = Q^{-1} - Q^{-1}[C_K^T \ D^T] \left( \begin{bmatrix} C_K \\ D \end{bmatrix} Q^{-1}[C_K^T \ D^T] \right)^{-1} \begin{bmatrix} C_K \\ D \end{bmatrix} Q^{-1}, K \in \mathcal{K}(x) \right\}.$$

We have the following proposition which states the first order sensitivity results associated with $s(\cdot)$. Its proof can be obtained through adapting the proofs in [18, Lemma 2.1] and [22, Theorem 1] with the help of change-of-variables.

**Proposition 2.** *For any $x \in \Re^n$, there exists a neighborhood $U$ of $x$ such that*
$$\mathcal{K}(y) \subseteq \mathcal{K}(x), \quad \mathcal{P}(y) \subseteq \mathcal{P}(x), \quad \forall y \in U.$$
*If $\mathcal{K}(y) \subseteq \mathcal{K}(x)$, it holds that $s(y) = s(x) + P(y - x), \forall P \in \mathcal{P}(y)$. Furthermore, let $I(x)$ be given in (10). Denote*
$$P_0 := Q^{-1} - Q^{-1}[C_{I(x)}^T \ D^T] \left( \begin{bmatrix} C_{I(x)} \\ D \end{bmatrix} Q^{-1}[C_{I(x)}^T \ D^T] \right)^{\dagger} \begin{bmatrix} C_{I(x)} \\ D \end{bmatrix} Q^{-1}.$$
*Then, $P_0 \in \mathcal{P}(x)$.*

Since $\mathcal{P}$ is obtained through generalizing the results of Han and Sun [18], we name it as "generalized HS-Jacobian". We end this section by showing that if the matrix $[C_K^T \ D^T]$ is a diagonal matrix with only 0-1 diagonal elements, then the procedure for computing a generalized HS-Jacobian $P \in \mathcal{P}(x)$ can be simplified greatly.

**Proposition 3.** *Let $\theta \in \Re^n$ be a given vector with each entry $\theta_i$ being 0 or 1 for all $i = 1, \ldots, n$. Let $\Theta = \operatorname{Diag}(\theta)$ and $\Sigma = I_n - \Theta$. It holds that*
$$P := Q^{-1} - Q^{-1} \Theta \left( \Theta Q^{-1} \Theta \right)^{\dagger} \Theta Q^{-1} = (\Sigma Q \Sigma)^{\dagger}. \tag{11}$$

*Proof.* We only consider the case that $\Theta \neq 0$ since the conclusion holds trivially if $\Theta = 0$. Define
$$\widehat{P} := I - Q^{-\frac{1}{2}} \Theta (\Theta Q^{-\frac{1}{2}} Q^{-\frac{1}{2}} \Theta)^{\dagger} \Theta Q^{-\frac{1}{2}}.$$
Then $P = Q^{-\frac{1}{2}} \widehat{P} Q^{-\frac{1}{2}}$. From [22, Lemma 1], we know that $\widehat{P} d = \Pi_{\operatorname{Null}(\Theta Q^{-\frac{1}{2}})}(d), \forall d \in \Re^n$. Since $\operatorname{Null}(\Theta Q^{-\frac{1}{2}}) = \operatorname{Ran}(Q^{\frac{1}{2}} \Sigma)$, we have that
$$\widehat{P} d = \Pi_{\operatorname{Ran}(Q^{\frac{1}{2}} \Sigma)}(d) = Q^{\frac{1}{2}} \Sigma (\Sigma Q \Sigma)^{\dagger} \Sigma Q^{\frac{1}{2}} d, \quad \forall d \in \Re^n.$$
Therefore, $P = \Sigma (\Sigma Q \Sigma)^{\dagger} \Sigma = (\Sigma Q \Sigma)^{\dagger}$, where the last equality follows from the fact that $\Sigma$ is a diagonal matrix with 0-1 diagonal elements. □

# 3 Efficient computations of the generalized Jacobian of $\text{Prox}_p(\cdot)$

In this section, we shall study the variational properties of the proximal mapping of the fused lasso regularizer $p$, namely the generalized Jacobian of $\text{Prox}_p$ and their efficient computations. Recall that the proximal mapping of $p$ is defined by

$$\text{Prox}_p(v) := \text{argmin}\left\{\lambda_1\|x\|_1 + \lambda_2\|Bx\|_1 + \frac{1}{2}\|x - v\|^2\right\}, \quad \forall v \in \Re^n,$$

where $\lambda_1, \lambda_2 \geq 0$ are given data. Denote also by $x_{\lambda_2}(v)$ the proximal mapping of $\lambda_2\|B \cdot \|_1$:

$$x_{\lambda_2}(v) := \text{argmin}\left\{\lambda_2\|Bx\|_1 + \frac{1}{2}\|x - v\|^2\right\}, \quad \forall v \in \Re^n. \tag{12}$$

Next we recall a key result in [14] concerning the computation of $\text{Prox}_p$.

**Proposition 4.** *[14, Proposition 1] Given $\lambda_1, \lambda_2 \geq 0$, it holds that*

$$\text{Prox}_p(v) = \text{Prox}_{\lambda_1\|\cdot\|_1}(x_{\lambda_2}(v)) = \text{sign}(x_{\lambda_2}(v)) \circ \max(|x_{\lambda_2}(v)| - \lambda_1, 0), \quad \forall v \in \Re^n.$$

The above proposition states that the proximal mapping of the fused lasso regularizer $\lambda_1\|\cdot\|_1 + \lambda_2\|B(\cdot)\|_1$ can be decomposed into the composition of the proximal mapping of $\lambda_1\|\cdot\|_1$ and the proximal mapping of $\lambda_2\|B(\cdot)\|_1$. See [44] for the extensions of the above result to other regularizers. Given $v \in \Re^n$, from Proposition 4, it is clear that the efficient computation of $\text{Prox}_p(v)$ mainly depends on the fast calculation of $x_{\lambda_2}(v)$. Fortunately, many efficient direct algorithms have been developed for the fast computation of $x_{\lambda_2}(v)$ [11, 9, 19]. Meanwhile, we note that the subgradient finding algorithm (SFA) designed in [24] is a fast iterative solver for computing $x_{\lambda_2}(v)$. The relative performance of most of the existing algorithms has been well documented in the recent paper [2], which appears to suggest that for large scale problems, the direct solver developed and implemented by Condat [9] has generally outperformed the other solvers. Hence, in our later numerical experiments, we will use Condat's algorithm and implementation[1] for computing $x_{\lambda_2}(v)$.

To study the variational properties of $\text{Prox}_p$, we first need the following lemma which provides an alternative way of computing $x_{\lambda_2}(\cdot)$ through the dual solution $z_{\lambda_2}(\cdot)$:

$$z_{\lambda_2}(u) := \text{argmin}\left\{\frac{1}{2}\|B^T z\|^2 - \langle z, u\rangle \mid \|z\|_\infty \leq \lambda_2\right\}, \quad \forall u \in \Re^{n-1}. \tag{13}$$

**Lemma 1.** *Given $\lambda_2 \geq 0$, it holds that $x_{\lambda_2}(v) = v - B^T z_{\lambda_2}(Bv), \forall v \in \Re^n$.*

*Proof.* The result follows directly from Fenchel's Duality Theorem [34, Theorem 31.3]. □

Given $v \in \Re^n$, by Proposition 4 and Lemma 1, we have that

$$\text{Prox}_p(v) = \text{Prox}_{\lambda_1\|\cdot\|_1}(x_{\lambda_2}(v)) = \text{Prox}_{\lambda_1\|\cdot\|_1}(v - B^T z_{\lambda_2}(Bv)). \tag{14}$$

Thus, if $z_{\lambda_2}(\cdot)$ is continuously differentiable near $Bv$ and $I - B^T z'_{\lambda_2}(Bv)B$ is nonsingular, then we would get by the chain-rule [40, Lemma 2.1] that

$$\partial_B \text{Prox}_p(v) = \left\{\Theta(I - B^T z'_{\lambda_2}(Bv)B) \mid \Theta \in \partial_B \text{Prox}_{\lambda_1\|\cdot\|_1}(x_{\lambda_2}(v))\right\},$$

---

[1] http://www.gipsa-lab.grenoble-inp.fr/~laurent.condat/download/condat_fast_tv.c

where $\partial_B$ denotes the B-subdifferential [8]. However, $z_{\lambda_2}(\cdot)$ may not be differentiable at $Bv$ and the above chain-rule is usually not available. Therefore, we need to define the generalized Jacobian of $\text{Prox}_p$ in a proper way. The technical details on how this can be done are presented next.

For any $u \in \Re^{n-1}$, since $BB^T$ is symmetric and positive definite, $z_{\lambda_2}(u)$ is the unique solution to the strongly convex QP (13). Therefore, the generalized HS-Jacobian of $z_{\lambda_2}$ can be obtained directly from the results developed in Section 2.2.

We start by defining some notation. Denote the active index set

$$I_z(v) := \{i \mid |(z_{\lambda_2}(Bv))_i| = \lambda_2,\ i = 1, \ldots, n-1\} \tag{15}$$

and a collection of index sets

$$\mathcal{K}_z(v) := \{\ K \subseteq \{1, \ldots, n-1\} \mid \text{supp}(Bx_{\lambda_2}(v)) \subseteq K \subseteq I_z(v)\}.$$

Note that from the optimality conditions for $z_{\lambda_2}(Bv)$, one can show that $\text{supp}(Bx_{\lambda_2}(v))$ is equal to the support of any optimal Lagrangian multiplier associated with the constraint $\|z\|_\infty \leq \lambda_2$. Define the mulifunction $\mathcal{P}_z : \Re^n \rightrightarrows \Re^{(n-1)\times(n-1)}$ by

$$\mathcal{P}_z(v) := \left\{\widehat{P} \in \Re^{(n-1)\times(n-1)} \mid \widehat{P} = (\Sigma_K BB^T \Sigma_K)^\dagger,\ K \in \mathcal{K}_z(v)\right\},$$

where $\Sigma_K = \text{Diag}(\sigma_K) \in \Re^{(n-1)\times(n-1)}$ with

$$(\sigma_K)_i = \begin{cases} 0, & \text{if}\quad i \in K, \\ 1, & \text{otherwise}, \end{cases}\quad i = 1, \ldots, n-1. \tag{16}$$

Note that according to Proposition 3, $\mathcal{P}_z(v)$ is exactly the generalized HS-Jacobian of $z_{\lambda_2}$ at $Bv$. Define the mulifunction $\mathcal{P}_x : \Re^n \rightrightarrows \Re^{n\times n}$ by

$$\mathcal{P}_x(v) := \left\{P \in \Re^{n\times n} \mid P = I - B^T \widehat{P} B,\ \widehat{P} \in \mathcal{P}_z(v)\right\}.$$

Here $\mathcal{P}_x(v)$ can be viewed as the generalized HS-Jacobian of $x_{\lambda_2}$ at $v$. More precisely, we can derive the following first order sensitivity results of $z_{\lambda_2}$ and $x_{\lambda_2}$.

**Proposition 5.** *For any $v \in \Re^n$, there exists a neighborhood $W$ of $v$ such that for all $w \in W$*

$$\mathcal{K}_z(w) \subseteq \mathcal{K}_z(v),\quad \mathcal{P}_z(w) \subseteq \mathcal{P}_z(v),\quad \mathcal{P}_x(w) \subseteq \mathcal{P}_x(v)$$

*and*

$$\begin{cases} z_{\lambda_2}(Bw) = z_{\lambda_2}(Bv) + \widehat{P}B(w-v),\quad \forall \widehat{P} \in \mathcal{P}_z(w), \\ x_{\lambda_2}(w) = x_{\lambda_2}(v) + P(w-v),\quad \forall P \in \mathcal{P}_x(w). \end{cases}$$

*Proof.* The desired results follow from Propositions 2 and 3, and [22, Lemma 1]. □

Next we show that we can derive an explicit formula to calculate the generalized Jacobian $P \in \mathcal{P}_x(v)$ when the special structure of $B$ is taken into consideration. For $2 \leq j \leq n$, define linear mappings $B_j : \Re^j \to \Re^{j-1}$ such that $B_j x = [x_1 - x_2; \ldots; x_{j-1} - x_j]$, $\forall x \in \Re^j$. With this notation, we can write $B = B_n$. The following lemma is needed for later discussions and can be proved through direct calculations.

**Lemma 2.** *For $2 \leq j \leq n$, it holds that*
$$T_j := I_j - B_j^T(B_jB_j^T)^{-1}B_j = \frac{1}{j}\mathbf{E}_j.$$

**Proposition 6.** *Let $\Sigma \in \Re^{(n-1)\times(n-1)}$ be an $N$-block diagonal matrix $\Sigma = \mathrm{Diag}(\Lambda_1, \ldots, \Lambda_N)$, where for $i = 1, \ldots, N$, $\Lambda_i$ is either the $n_i$ by $n_i$ zero matrix $\mathbf{O}_{n_i}$ or the $n_i$ by $n_i$ identity matrix $I_{n_i}$ and any two consecutive blocks cannot be of the same type. Denote $J := \{j \mid \Lambda_j = I_{n_j}, j = 1, \ldots, N\}$. Then, it holds that*
$$\Gamma := I_n - B^T(\Sigma BB^T\Sigma)^\dagger B = \mathrm{Diag}(\Gamma_1, \ldots, \Gamma_N),$$
*where for $i = 1, \ldots, N$,*
$$\Gamma_i = \begin{cases} \dfrac{1}{n_i+1}\mathbf{E}_{n_i+1}, & \text{if } i \in J, \\ I_{n_i}, & \text{if } i \notin J \text{ and } i = 1 \text{ or } N, \\ I_{n_i-1}, & \text{otherwise} \end{cases} \tag{17}$$
*with the convention $I_0 = \emptyset$. Moreover, $\Gamma = H + UU^T = H + U_J U_J^T$, where $H \in \Re^{n \times n}$ is an $N$-block diagonal matrix given by $H = \mathrm{Diag}(\Upsilon_1, \ldots, \Upsilon_N)$ with*
$$\Upsilon_i = \begin{cases} \mathbf{O}_{n_i+1}, & \text{if } i \in J, \\ I_{n_i}, & \text{if } i \notin J \text{ and } i = 1 \text{ or } N, \\ I_{n_i-1}, & \text{otherwise} \end{cases}$$
*and $U \in \Re^{n \times N}$ with its $(k,j)$-th entry given by*
$$U_{k,j} = \begin{cases} \dfrac{1}{\sqrt{n_j+1}}, & \text{if } \sum_{t=1}^{j-1} n_t + 1 \leq k \leq \sum_{t=1}^{j} n_t + 1, \text{ and } j \in J, \\ 0, & \text{otherwise} \end{cases} \tag{18}$$
*and $U_J$ consists of the nonzero columns of $U$, i.e., the columns whose indices are in $J$.*

*Proof.* Note that $(\Sigma BB^T\Sigma)^\dagger = \mathrm{Diag}(T_1, \ldots, T_N)$, where for $i = 1, \ldots, N$,
$$T_i = \begin{cases} (B_{n_i}B_{n_i}^T)^{-1} & \text{if } \Lambda_i = I_{n_i}, \\ \mathbf{O}_{n_i}, & \text{otherwise.} \end{cases}$$
Then by Lemma 2 and the structure of $B$, we can obtain the desired results through some direct calculations. □

Define the multifunction $\mathcal{M} : \Re^n \rightrightarrows \Re^{n \times n}$ by
$$\mathcal{M}(v) := \left\{ M \in \Re^{n \times n} \mid M = \Theta P,\ \Theta \in \partial_B \mathrm{Prox}_{\lambda_1 \|\cdot\|_1}(x_{\lambda_2}(v)),\ P \in \mathcal{P}_x(v) \right\}. \tag{19}$$

Let $v \in \Re^n$ be an arbitrary point. The set $\mathcal{M}(v)$ is exactly the generalized Jacobian of $\mathrm{Prox}_p$ at $v$ to be used in this paper. In numerical implementations, one needs to construct at least one element in $\mathcal{M}(v)$ explicitly. This can be done in the following manner. Firstly, denote $\Theta = \mathrm{Diag}(\theta)$ with
$$\theta_i = \begin{cases} 0, & \text{if } |(x_{\lambda_2}(v))_i| \leq \lambda_1, \\ 1, & \text{otherwise}, \quad i = 1, \ldots, n. \end{cases} \tag{20}$$

Then, let $I_z(v)$ be given as in (15) and $P = I_n - B^T(\Sigma BB^T\Sigma)^\dagger B$, where $\Sigma = \text{Diag}(\sigma) \in \Re^{(n-1)\times(n-1)}$ with

$$\sigma_i = \begin{cases} 0, & \text{if } i \in I_z(v), \\ 1, & \text{otherwise}, \quad i = 1, \ldots, n-1. \end{cases} \tag{21}$$

Obviously, $\Theta \in \partial_B \text{Prox}_{\lambda_1 \|\cdot\|_1}(x_{\lambda_2}(v))$ and $P \in \mathcal{P}_x(v)$. Therefore,

$$M := \Theta P \in \mathcal{M}(v). \tag{22}$$

The following main theorem of this section shows why $\mathcal{M}(v)$ can be indeed regarded as the generalized Jacobian of $\text{Prox}_p$ at $v$.

**Theorem 1.** *Let $\lambda_1, \lambda_2 \geq 0$ be nonnegative numbers and $v \in \Re^n$. Then, $\mathcal{M}$ is a nonempty and compact valued and upper-semicontinuous multifunction and for any $M \in \mathcal{M}(v)$, $M$ is symmetric and positive semidefinite. Moreover, there exists a neighborhood $W$ of $v$ such that for all $w \in W$*

$$\text{Prox}_p(w) - \text{Prox}_p(v) - M(w - v) = 0, \quad \forall M \in \mathcal{M}(w). \tag{23}$$

*Proof.* It is obvious that the point-to-set mapping $\mathcal{M}$ has nonempty compact images. The upper semicontinuity of $\mathcal{M}$ follows from the Lipschitz continuity of $x_{\lambda_2}()$ and the upper semicontinuity of the B-subdifferential mapping $\partial_B \text{Prox}_{\lambda_1 \|\cdot\|_1}(\cdot)$ and the inclusion property on $\mathcal{P}_x(\cdot)$ obtained in Proposition 5. Since $\text{Prox}_{\lambda_1 \|\cdot\|_1}(\cdot)$ is piecewise affine and $x_{\lambda_2}(\cdot)$ is Lipschitz continuous, equation (23) follows easily from Proposition 5 and [12, Theorem 7.5.17]. Thus, we only need to prove that for any $v \in \Re^n$ and $M \in \mathcal{M}(v)$, $M \in \mathcal{S}_+^n$, the set of $n \times n$ symmetric and positive semidefinite matrices. Indeed, for any $M \in \mathcal{M}(v)$, there exist $\Theta \in \partial_B \text{Prox}_{\lambda_1 \|\cdot\|}(x_{\lambda_2}(v))$ and $K \in \mathcal{K}_z(v)$ such that $M = \Theta(I - B^T(\Sigma_K BB^T\Sigma_K)^\dagger B)$ with $\Sigma_K$ given in (16). From Proposition 6, we have $I - B^T(\Sigma_K BB^T\Sigma_K)^\dagger B = \text{Diag}(\Gamma_1, \ldots, \Gamma_N)$ with $\Gamma_i$ given in (17). Note that $\Theta$ can also be decomposed with the same pattern as $\Gamma$, i.e.,

$$\Theta = \text{Diag}(\Theta_1, \ldots, \Theta_N).$$

Thus $M = \text{Diag}(\Theta_1 \Gamma_1, \ldots, \Theta_N \Gamma_N)$. Let $J := \{j \mid \Gamma_j \text{ is not an identity matrix}, 1 \leq j \leq N\}$. Since $\text{supp}(Bx_{\lambda_2}(v)) \subseteq K$, we have that

$$\Theta_j = \mathbf{O}_{n_j+1} \quad \text{or} \quad I_{n_j+1}, \quad \forall j \in J,$$

which implies $\Theta_j \Gamma_j \in \mathcal{S}_+^{n_j+1}, \forall j \in J$. Therefore, $M \in \mathcal{S}_+^n$ and the proof is completed. □

Theorem 1 indicates that for an arbitrary constant $\gamma > 0$, the function $\text{Prox}_p$ is in fact $\gamma$-order semismooth on $\Re^n$ with respect to the multifunction $\mathcal{M}$ in the sense of the following definition of semismoothness from [26, 20, 32, 39].

**Definition 1** (Semismoothness). *Let $\mathcal{O} \subseteq \Re^n$ be an open set, $\mathcal{K} : \mathcal{O} \subseteq \Re^n \rightrightarrows \Re^{m \times n}$ be a nonempty and compact valued, upper-semicontinous set-valued mapping and $F : \mathcal{O} \to \Re^m$ be a locally Lipschitz continuous function. $F$ is said to be semismooth at $x \in \mathcal{O}$ with respect to the multifunction $\mathcal{K}$ if $F$ is directionally differentiable at $x$ and for any $V \in \mathcal{K}(x + \Delta x)$ with $\Delta x \to 0$,*

$$F(x + \Delta x) - F(x) - V\Delta x = o(\|\Delta x\|).$$

*Let $\gamma$ be a positive constant. $F$ is said to be $\gamma$-order (strongly, if $\gamma = 1$) semismooth at $x$ with respect to $\mathcal{K}$ if $F$ is directionally differentiable at $x$ and for any $V \in \mathcal{K}(x + \Delta x)$ with $\Delta x \to 0$,*

$$F(x + \Delta x) - F(x) - V \Delta x = O(\|\Delta x\|^{1+\gamma}).$$

*$F$ is said to be a semismooth (respectively, $\gamma$-order semismooth, strongly semismooth) function on $\mathcal{O}$ with respect to $\mathcal{K}$ if it is semismooth (respectively, $\gamma$-order semismooth, strongly semismooth) everywhere in $\mathcal{O}$ with respect to $\mathcal{K}$.*

**Remark 2.** *Note that as a Lipschitz continuous piecewise affine function, $\text{Prox}_p$ is $\gamma$-order semismooth on $\Re^n$ with respect to the Clarke generalized Jacobian $\partial \text{Prox}_p$ for any given $\gamma > 0$ [8].*

## 4 A semismooth Newton based augmented Lagrangian method for fused lasso regularized least squares problems

In this section, we present the backbone of the level-set method – a highly efficient algorithm for solving the fused lasso regularized least squares problems arising at each iteration of the method. Critical numerical issues concerning its efficient implementations will also be discussed.

Given $A \in \Re^{m \times n}$, $b \in \Re^m$, $\lambda_1, \lambda_2 \geq 0$, note that the subproblems of the level-set method can be written as follows

$$(\mathbf{P}) \quad \max \left\{ f(x) := -\frac{1}{2} \|Ax - b\|^2 - p(x) \right\},$$

where the fused lasso regularizer $p$ is given by $p(x) = \lambda_1 \|x\|_1 + \lambda_2 \|Bx\|_1$, $\forall x \in \Re^n$. It can shown readily that the dual of ($\mathbf{P}$) is given by

$$(\mathbf{D}) \quad \min \left\{ g(y) := \frac{1}{2} \|y\|^2 + \langle y, b \rangle + p^*(-A^T y) \right\}.$$

Now, we derive the augmented Lagrangian function for the unconstrained minimization problem ($\mathbf{D}$) following the framework presented in [37, Examples 11.46 and 11.57]. Firstly, we identify problem ($\mathbf{D}$) with the problem of minimizing $g(y) = \tilde{g}(y, 0)$ over $\Re^m$ for

$$\tilde{g}(y, \xi) := h^*(y) + p^*(-A^T y + \xi), \quad \forall (y, \xi) \in \Re^m \times \Re^n.$$

Note that $\tilde{g}$ is jointly convex in $(y, \xi)$. The Lagrangian function $l : \Re^m \times \Re^n \to [-\infty, +\infty]$ associated with ($\mathbf{D}$) is given by

$$l(y; x) = \inf_\xi \{\tilde{g}(y, \xi) - \langle x, \xi \rangle\} = \frac{1}{2} \|y\|^2 + \langle y, b \rangle - \langle A^T y, x \rangle - p(x). \tag{24}$$

Given $\sigma > 0$, the augmented Lagrangian function corresponding to ($\mathbf{D}$) can be obtained by

$$\mathcal{L}_\sigma(y; x) = \inf_\xi \left\{ \tilde{g}(y, \xi) - \langle x, \xi \rangle + \frac{\sigma}{2} \|\xi\|^2 \right\} = \sup_s \left\{ l(y; s) - \frac{1}{2\sigma} \|s - x\|^2 \right\}$$

$$= \frac{1}{2} \|y\|^2 + \langle y, b \rangle + \inf_s \left\{ p^*(s) - \langle x, A^T y + s \rangle + \frac{\sigma}{2} \|A^T y + s\|^2 \right\}, \quad \forall (y, x) \in \Re^m \times \Re^n.$$

Now, we can present the detailed steps of algorithm SSNAL for solving (**D**) as follows. The algorithm is termed as SSNAL since a semismooth Newton method will be employed to solve the subproblems involved in the inexact augmented Lagrangian method [35].

---

**Algorithm** SSNAL: **A semismooth Newton based augmented Lagrangian method for (D).**

Let $\sigma_0 > 0$ be a given parameter. Choose $(y^0, x^0) \in \Re^m \times \Re^n$. For $k = 0, 1, \ldots$, perform the following steps in each iteration:

**Step 1.** Compute
$$y^{k+1} \approx \arg\min\{\Psi_k(y) := \mathcal{L}_{\sigma_k}(y; x^k)\}. \tag{25}$$

**Step 2.** Compute $x^{k+1} = \text{Prox}_{\sigma_k p}(x^k - \sigma_k A^T y^{k+1})$.

**Step 3.** Update $\sigma_{k+1} \uparrow \sigma_\infty \leq \infty$ .

---

Since $\Psi_k$ is strongly convex and differentiable, the standard stopping criteria studied in [36, 35] for approximately solving (25) can be simplified to

(A)    $\|\nabla \Psi_k(y^{k+1})\| \leq \varepsilon_k/\sqrt{\sigma_k}, \quad \sum_{k=0}^{\infty} \varepsilon_k < \infty,$

(B1)   $\|\nabla \Psi_k(y^{k+1})\| \leq (\delta_k/\sqrt{\sigma_k}) \|\text{Prox}_{\sigma_k p}(x^k - \sigma_k A^T y^{k+1}) - x^k\|, \quad \sum_{k=0}^{\infty} \delta_k < +\infty,$

(B2)   $\|\nabla \Psi_k(y^{k+1})\| \leq (\delta'_k/\sigma_k) \|\text{Prox}_{\sigma_k p}(x^k - \sigma_k A^T y^{k+1}) - x^k\|, \quad 0 \leq \delta'_k \to 0.$

The global and local convergence results of the above algorithm have been extensively studied in [36, 35, 25, 21]. We also refer the reader to the recent paper [10] by Cui et al. for more discussions on new implementable stopping criteria and the superlinear convergence of the ALM. Here we will list some of them without proofs.

**Theorem 2.** *Suppose that the solution set to* (**P**) *is nonempty. Let* $\{(y^k, z^k, x^k)\}$ *be the infinite sequence generated by Algorithm* SSNAL *with stopping criterion* (A). *Then, the sequence* $\{x^k\}$ *is bounded and converges to an optimal solution of* (**P**). *Moreover,* $\{y^k\}$ *is also bounded and converges to the unique optimal solution* $y^*$ *of* (**D**).

*If Algorithm* SSNAL *is executed under stopping criteria* (A) *and* (B1), *then* $\{x^k\}$ *converges asymptotically Q-superlinearly. If in addition, the stopping criterion* (B2) *is also used, then* $\{y^k\}$ *converges asymptotically R-superlinearly.*

## 4.1 Solving the subproblems in SSNAL

As the reader may observe, the most expensive computation in each iteration of the augmented Lagrangian method is to solve the subproblem in Step 1. Let $\sigma > 0$ and $\tilde{x} \in \mathcal{X}$ be fixed. We shall propose an efficient semismooth Newton algorithm to solve the following inner problem involved in each iteration of the above SSNAL method (25):

$$\min_{y \in \Re^m} \Psi(y) := \mathcal{L}_\sigma(y; \tilde{x}). \tag{26}$$

Since $\Psi$ is a strongly convex function on $\Re^m$, we have that, for any $\alpha \in \Re$, the level set $\mathcal{L}_\alpha := \{y \in \Re^m \mid \Psi(y) \leq \alpha\}$ is a closed and bounded convex set. Moreover, problem (25) admits a unique optimal solution denoted as $\bar{y}$.

For $y \in \Re^m$, denote $x(y) := \tilde{x} - \sigma A^T y$. Note that $\Psi$ is continuously differentiable on $\Re^m$ with

$$\nabla \Psi(y) = y - A\operatorname{Prox}_{\sigma p}(x(y)), \quad \forall y \in \Re^m.$$

Thus, the unique optimal solution $\bar{y}$ of (26) can be obtained via solving the following nonsmooth equation

$$\nabla \Psi(y) = 0. \tag{27}$$

We can define the following multifunction $\mathcal{V}: \Re^m \rightrightarrows \Re^{m \times m}$ by

$$\mathcal{V}(y) := \{V \in \Re^{m \times m} \mid V = I_m + \sigma A M A^T, M \in \mathcal{M}(x(y))\},$$

where $\mathcal{M}$ is the multifunction defined in (19). In contrast to the cases studied in [46, 21], here the calculation of the Clarke generalized Jacobian of $\partial \operatorname{Prox}_{\sigma p}$ [8] is much more involved given the inherited composition structure of the nonsmooth fused lasso regularizer $p$. Thus we propose to use $\mathcal{M}$ to replace $\partial \operatorname{Prox}_{\sigma p}$. Obviously, the multifunction $\mathcal{V}$ is also nonempty and compact valued and upper-semicontinous and $\nabla \Psi$ is $\gamma$-order semismooth on $\Re^m$ with respect to $\mathcal{V}$ for any $\gamma > 0$. Moreover, from Theorem 1, we know that every element in $\mathcal{V}(y)$ is symmetric and positive definite.

We apply the following semismooth Newton (SSN) method to solve (27), and expect to get a fast superlinear or even quadratic convergence.

---

**Algorithm** SSN: **A semismooth Newton algorithm for solving (27)** (SSN($y^0, \tilde{x}, \sigma$)).

Given $\mu \in (0, 1/2)$, $\bar{\eta} \in (0,1)$, $\tau \in (0,1]$, and $\delta \in (0,1)$, choose $y^0 \in \Re^m$. Iterate the following steps for $j = 0, 1, \ldots$.

**Step 1.** Let $M_j$ be a generalized Jacobian of $p$ at $\tilde{x} - \sigma A^T y^j$ as defined in (22). Set $V_j := I_m + \sigma A M_j A^T$. Solve the following linear system

$$V_j d = -\nabla \Psi(y^j) \tag{28}$$

exactly or by the conjugate gradient (CG) algorithm to find $d^j$ such that $\|V_j d^j + \nabla \Psi(y^j)\| \leq \min(\bar{\eta}, \|\nabla \Psi(y^j)\|^{1+\tau})$.

**Step 2.** (Line search) Set $\alpha_j = \delta^{m_j}$, where $m_j$ is the first nonnegative integer $m$ for which

$$\Psi(y^j + \delta^m d^j) \leq \Psi(y^j) + \mu \delta^m \langle \nabla \Psi(y^j), d^j \rangle.$$

**Step 3.** Set $y^{j+1} = y^j + \alpha_j d^j$.

---

In order to study the convergence results for the above SSN algorithm, the following proposition will be needed in the sequel. It is in fact an extension of [29, Theorem 2.1].

**Proposition 7.** *Let $\theta : \Omega \to \Re$, with $\Omega$ open, be a continuously differentiable function and its gradient $\nabla \theta : \Omega \to \Re^n$ is locally Lipschitz in $\Omega$. If $\nabla \theta$ is semismooth at a point $x \in \Omega$ with respect to a nonempty and compact valued and upper-semicontinuous multifunction $\mathcal{K} : \Omega \rightrightarrows \mathcal{S}^n$, then for any $V \in \mathcal{K}(x + d)$ with $d \to 0$, we have*

$$\theta(x + d) - \theta(x) - \langle \nabla \theta(x), d \rangle - \frac{1}{2} \langle d, V d \rangle = o(\|d\|^2).$$

*Proof.* From [38] and the semismoothness of $\nabla\theta$ at $x$, we know that $Vd-(\nabla\theta)'(x;d) = o(\|d\|)$, $\forall d \to 0$ and $V \in \mathcal{K}(x+d)$. Then we get the desired limit by following the proof of [29, Theorem 2.1]. □

**Theorem 3.** *Let $\{y^j\}$ be the infinite sequence generated by Algorithm* SSN. *Then $\{y^j\}$ converges to the unique optimal solution $\bar{y}$ of problem* (26) *and $\|y^{j+1} - \bar{y}\| = O(\|y^j - \bar{y}\|^{1+\tau})$.*

*Proof.* Since by [46, Proposition 3.3], $d^j$ is always a descent direction, Algorithm SSN is well-defined. The strong convexity of $\Psi(\cdot)$ and [46, Theorem 3.4] imply that $\{y^j\}$ converges to the unique optimal solution $\bar{y}$ of problem (26).

Note that every $V \in \mathcal{V}(\bar{y})$ is symmetric and positive definite. Since $\mathcal{V}$ is upper-semicontinuous, from [12, Lemma 7.5.2], we have that for all $j$ sufficiently large, $\{\|V_j^{-1}\|\}$ is uniformly bounded. Since $\nabla\Psi$ is strongly semismooth with respect to $\mathcal{V}$, similar to the proof for [46, Theorem 3.5], it can be shown that for all $j$ sufficiently large,

$$\|y^j + d^j - \bar{y}\| = \mathcal{O}(\|y^j - \bar{y}\|^{1+\tau}), \tag{29}$$

and for some constant $\hat{\delta} > 0$, $-\langle\nabla\Psi(y^j), d^j\rangle \geq \hat{\delta}\|d^j\|^2$. Based on (29), Proposition 7 and [12, Proposition 8.3.18], we know that for $\mu \in (0, 1/2)$, there exists an integer $j_0$ such that for all $j \geq j_0$,

$$\Psi(y^j + d^j) \leq \Psi(y^j) + \mu\langle\nabla\Psi(y^j), d^j\rangle,$$

i.e., $y^{j+1} = y^j + d^j$ for all $j \geq j_0$. This, together with (29), completes the proof. □

## 4.2 Efficient implementations of the semismooth Newton method by exploiting second-order sparsity

When Algorithm SSN is used to solve the subproblem (26), the most expensive step is the computation of the search direction $d^j$ from the linear system (28). Given $(x, y) \in \Re^n \times \Re^m$ and $\sigma > 0$, the Newton system (28) associated with the fused lasso problem is given by

$$(I_m + \sigma AMA^T)d = -\nabla\Psi(y), \tag{30}$$

where $M \in \mathcal{M}(\tilde{x})$ is given in (22) with $\tilde{x} := x - \sigma A^T y$. We note that in the case of the standard lasso problem, the counterpart of $M$ is a diagonal matrix but here the structure of $M$ is much more complex. Since $M$ is an $n \times n$ matrix, the costs of naively computing $AMA^T$ and the matrix-vector multiplication $AMA^T d$ for a given vector $d \in \Re^m$ are $\mathcal{O}(n^2 m + nm^2)$ and $\mathcal{O}(n^2 + mn)$, respectively. Thus, neither the Cholesky factorization nor the conjugate gradient method would be efficient for solving the above linear system when $n$ and/or $m$ is large. In fact, in the high-dimensional setting where the number of features $n$ is usually of the order $10^5$, it would be impossible to store $M$ in the RAM of a standard desktop computer when $M$ is dense.

Fortunately, in the previous section, we have already extracted and analysed the structure of $M$. From equation (22), we know that

$$M = \Theta P \quad \text{with} \quad P = I_n - B^T(\Sigma BB^T\Sigma)^\dagger B,$$

where diagonal matrices $\Theta$, $\Sigma$ are given in (20) and (21), respectively. Note that $\Sigma$ is an $N$-block diagonal matrix $\Sigma = \text{Diag}(\Lambda_1, \ldots, \Lambda_N)$ with each $\Lambda_i \in \Re^{n_i \times n_i}$ being either the zero matrix or the identity matrix with different types for any two consecutive blocks. Let $J := \{j \mid \Lambda_j = I_{n_j}, j = $ $1, \ldots, N\}$. It can be seen from Proposition 6 that $P$ could be expressed as the sum of a low rank matrix and a diagonal matrix. More specifically, we have that

$$P = H + UU^T = H + U_J U_J^T,$$

where $H \in \Re^{n \times n}$ is an $N$-block diagonal matrix given by

$$H = \text{Diag}(\Upsilon_1, \ldots, \Upsilon_N) \quad \text{with} \quad \Upsilon_i = \begin{cases} \mathbf{O}_{n_i+1}, & \text{if } i \in J, \\ I_{n_i}, & \text{if } i \notin J \text{ and } i = 1 \text{ or } N, \\ I_{n_i-1}, & \text{otherwise} \end{cases}$$

and $U \in \Re^{n \times N}$ with its $(k, j)$-th entry given by

$$U_{k,j} = \begin{cases} \dfrac{1}{\sqrt{n_j + 1}}, & \text{if } \sum_{t=1}^{j-1} n_t + 1 \leq k \leq \sum_{t=1}^{j} n_t + 1, \text{ and } j \in J, \\ 0, & \text{otherwise} \end{cases}$$

and $U_J$ consists of the nonzero columns of $U$, i.e., the columns whose indices are in $J$.

Since $M$ is a symmetric matrix and $\Theta$, $H$ are diagonal matrices with only 0-1 diagonal elements, it holds that

$$M = \Theta(H + U_J U_J^T) = \Theta(H + U_J U_J)^T \Theta, \quad \Theta^2 = \Theta, \quad H^2 = H, \quad \Theta H = \Theta H \Theta.$$

Therefore,

$$A\Theta H A^T = A\Theta H \Theta A^T, \quad A\Theta(U_J U_J^T) A^T = A\Theta(U_J U_J^T) \Theta A^T.$$

Define the following index sets

$$\alpha := \{i \mid \theta_i = 1, i \in \{1, \ldots, n\}\}, \quad \beta := \{i \mid h_i = 0, i \in \alpha\},$$

where $\theta_i$ and $h_i$ are the $i$-th diagonal entries of matrices $\Theta$ and $H$, respectively. Then, we have

$$A\Theta H A^T = A\Theta H \Theta A^T = A_\alpha H A_\alpha^T = A_\beta A_\beta^T,$$

where $A_\alpha \in \Re^{m \times |\alpha|}$ and $A_\beta \in \Re^{m \times |\beta|}$ are two sub-matrices obtained from $A$ by extracting those columns with indices in $\alpha$ and $\beta$, respectively. Meanwhile, we have

$$A\Theta(U_J U_J^T) A^T = A\Theta(U_J U_J^T) \Theta A^T = A_\alpha \widetilde{U} \widetilde{U}^T A_\alpha^T,$$

where $\widetilde{U} \in \Re^{|\alpha| \times r}$ is a sub-matrix obtained from $\Theta U_J$ by extracting those rows with indices in $\alpha$ and the zeros columns in $\Theta U_J$ are removed. Here $r$ is the number of columns of $\widetilde{U}$. Note that $|\beta| \leq |\alpha| \leq n$ and $r \leq |J| \leq n$. Therefore, by exploiting the structure in $M$, we can express $AMA^T$ in the following form

$$AMA^T = A_\beta A_\beta^T + A_\alpha \widetilde{U} \widetilde{U}^T A_\alpha^T. \tag{31}$$

We call the above structure of $AMA^T$ and that of $I_m + \sigma AMA^T$ inherited from $M$ as second-order structured sparsity. The term "second-order" is used because $I_m + \sigma AMA^T$ can be viewed as a generalized Hessian of $\Psi$ at the given point $y$.

From the structure uncovered in (31), we can reduce the costs of computing $AMA^T$ and $AMA^T d$ for a given vector $d$ to $\mathcal{O}(m|\alpha|(m+r))$ and $\mathcal{O}((m+r)|\alpha|)$, respectively. Due to the presence of the fused lasso regularizer, $|\alpha|$, $|\beta|$ and $r$ usually are much smaller than $n$. Thus, by carefully exploring the special "low rank + diagonal" structure and the hidden sparsity in $M$, we are able to achieve a significant reduction in the cost of solving the linear system (30). More specifically, the total computational cost of using the Cholesky factorization to solve the linear system is reduced from $\mathcal{O}(m^3) + \mathcal{O}(m^2 n)$ to $\mathcal{O}(m^3) + \mathcal{O}(m^2|\alpha|(1+r/m))$. We note that a similar but much simpler reduction has been exploited in [21] for solving the classical lasso problems. See Figure 1 in [21, Section 3.3] for a graphical illustration.

The computational cost for solving (30) can be further reduced if the Sherman-Morrison-Woodbury formula [17] is properly used. When $m$ is large and the optimal solution is indeed sparse, it is very likely to have $r + |\beta| \ll m$ and/or $|\beta| + |\alpha| \ll m$. If $r + |\beta| \ll m$, let $W := [A_\beta, A_\alpha \widetilde{U}] \in \Re^{m \times (r+|\beta|)}$. Then, we know that

$$AMA^T = WW^T \quad \text{and} \quad (I_m + \sigma AMA^T)^{-1} = I_m - W(\sigma^{-1} I_{r+|\beta|} + W^T W)^{-1} W^T.$$

That is, instead of factorizing an $m \times m$ matrix, we only need to factorize an $(r+|\beta|) \times (r+|\beta|)$ matrix. In this case, the computational cost is reduced from $\mathcal{O}(m^3) + \mathcal{O}(m^2|\alpha|(1+r/m))$ to $\mathcal{O}((r+|\beta|)^3) + \mathcal{O}(m(r+|\beta|)^2)$. Similarly, if $|\beta| + |\alpha| \ll m$, we have the following decomposition

$$AMA^T = W_1 W_2^T \text{ with } W_1 := [A_\beta, A_\alpha \widetilde{U}\widetilde{U}^T] \in \Re^{m \times (|\alpha|+|\beta|)}, \ W_2 := [A_\beta, A_\alpha] \in \Re^{m \times (|\alpha|+|\beta|)}.$$

Using the above decomposition of $AMA^T$, we get

$$(I_m + \sigma AMA^T)^{-1} = I_m - W_1(\sigma^{-1} I_{|\alpha|+|\beta|} + W_2^T W_1)^{-1} W_2^T.$$

Thus, we only need to factorize an $(|\alpha|+|\beta|) \times (|\alpha|+|\beta|)$ matrix and the total computational cost is merely $\mathcal{O}((|\alpha|+|\beta|)^3) + \mathcal{O}(m(|\alpha|+|\beta|)^2)$ instead of the original $\mathcal{O}(m^3 + m^2|\alpha| + mr|\alpha|)$. In either way, we are able to greatly reduce the computational cost for solving the linear system (30).

## 5 Numerical experiments

In this section, we first evaluate the performance of SSNAL for solving fused lasso regularized least squares problems. Next, we demonstrate the power of SSNAL in solving regularized least squares subproblems within the level-set method for solving the least squares constrained fused lasso problems (1). All our computational results are obtained by running MATLAB (version 9.0) on a windows workstation (12-core, Intel Xeon E5-2680 @ 2.50GHz, 128 G RAM).

### 5.1 Numerical results for fused lasso regularized least squares problems

For solving the fused lasso regularized least squares problems (2), we will compare SSNAL with the state-of-the-art solver for fused lasso problems SLEP[2] [23] (which is based on Nesterov's accelerated proximal gradient method [27, 3]) and the popular alternating direction method of multipliers (ADMM) [15, 16]. For comparison purposes, we also test the inexact ADMM (iADMM) [7] and the linearized ADMM (LADMM) [45]. As is already mentioned in Section 3, the efficiency of computing the proximal mapping $\text{Prox}_p$ depends critically on the availability of a fast solver for computing the proximal mapping $\text{Prox}_{\lambda_2\|B(\cdot)\|_1}$ of the TV-norm $\|B(\cdot)\|_1$. From the results presented in [2], it appears that currently the most efficient code for computing $\text{Prox}_{\lambda_2\|B(\cdot)\|_1}$ is Condat's direct algorithm[3]. Hence, we use this direct algorithm in all the tested algorithms in the computation of $\text{Prox}_p$. In particular, we enhanced the performance of SLEP by replacing the subgradient finding algorithm for computing $\text{Prox}_{\lambda_2\|B(\cdot)\|_1}$ in SLEP by Condat's algorithm. While SLEP is used to sovle the primal problem (**P**), the ADMM type of methods are used to solve the following variants of problem (**D**)

$$\min\left\{\frac{1}{2}\|y\|^2 + \langle y,\,b\rangle + p^*(u) \,|\, A^T y + u = 0\right\}.$$

---

[2] http://www.yelab.net/software/SLEP/

The main difference between ADMM and iADMM is how the linear systems corresponding to $y$ in the subproblems are solved. While ADMM solves the linear systems exactly (up to machine precision) using a direct method, iADMM can use an iterative solver such as the preconditioned conjugate gradient (PCG) method to solve the linear systems inexactly. Clearly, when the linear system is large and using a direct solver is expensive, iADMM would be preferred over ADMM. We have implemented SSNAL, ADMM, iADMM and LADMM in MATLAB. In particular, in our implementation of ADMM type methods, we set the step-length to be 1.618.

In our numerical experiments, the regularization parameters $\lambda_1$ and $\lambda_2$ in the fused lasso problem (**P**) are chosen as

$$\lambda_1 = \alpha_1 \|A^T b\|_\infty \quad \text{and} \quad \lambda_2 = \alpha_2 \lambda_1,$$

where $0 < \alpha_1 < 1$ and $\alpha_2 > 0$. We measure the accuracy of an approximate optimal solution $\tilde{x}$ for (**P**) by using the following relative KKT residual:

$$\eta = \frac{\|\tilde{x} - \text{Prox}_p(\tilde{x} - A^T(A\tilde{x} - b))\|}{1 + \|\tilde{x}\| + \|A\tilde{x} - b\|}.$$

We stop the tested algorithms when $\eta \leq 10^{-6}$. The algorithms will also be stopped when they reach the pre-set maximum number of iterations (100 for SSNAL, and 20,000 for SLEP, ADMM, iADMM and LADMM) or the maximum computation time of 7 hours. All the parameters for SLEP are set to the default values.

### 5.1.1 Numerical results for high-dimentional biomedical datasets

In this subsection, we compare the algorithms on the test instances $(A, b)$ obtained from Kent Ridge Biomedical Data Set Repository[4]. During the data collection process, we normalize the matrix $A$ to have columns with at most unit norm. We extract 10 instances, namely `DLBCLH`, `DLBCN`, `DLBCLS`, `lungH1`, `lungH2`, `lungM`, `lungO`, `NervousSystem`, `ovarianP` and `overianS`. All the instances are in the high-dimension-low-sample size setting.

We choose the parameters $\alpha_1 \in \{10^{-3}, 10^{-4}\}$ and $\alpha_2 \in \{10, 2, 0.02, 0.01\}$, respectively. That is, we test 80 instances in total. In our test, ADMM and iADMM are able to successfully solve 68 and 73 instances, respectively, while SLEP fails to solve any of the instances to the required accuracy of $10^{-6}$. SSNAL is the only algorithm that can solve all 80 instances successfully. We note that the poor performance of SLEP and LADMM is closely related to the fact that for all these examples, the Lipschitz constants $\|A\|_2^2$ for the quadratic functions $\frac{1}{2}\|Ax - b\|^2$ are all rather large.

---

[3] http://www.gipsa-lab.grenoble-inp.fr/~laurent.condat/download/condat_fast_tv.c
[4] http://leo.ugr.es/elvira/DBCRepository/index.html

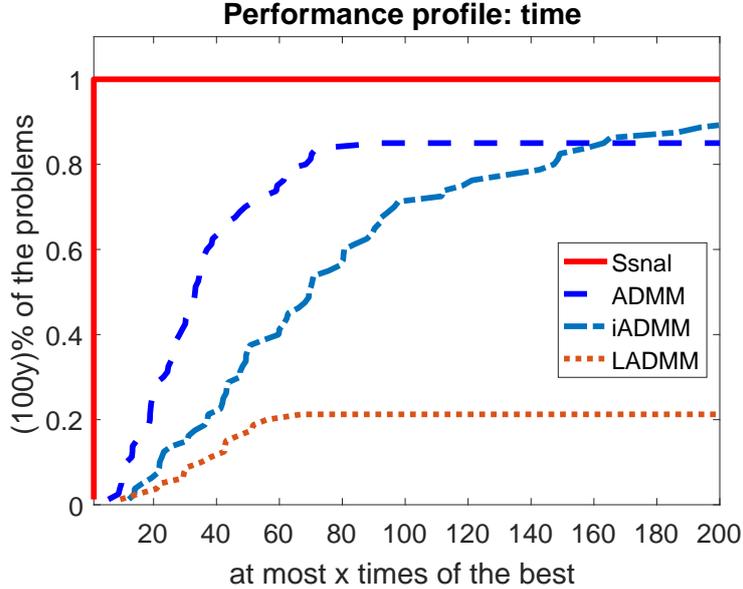

Figure 1: Performance profiles of S<small>SNCG</small>, ADMM, iADMM and LADMM on biomedical datasets.

Table 1 reports the detailed numerical results for S<small>SNAL</small>, SLEP, ADMM, iADMM and LADMM in solving some of the larger instances in the biomedical datasets. (The full results for this and subsequent tables are available at http://www.math.nus.edu.sg/~mattohkc/papers/fusedlassotables.pdf.) In the table, $m$ denotes the number of samples, $n$ denotes the number of features, "nnz($x$)" and "nnz($Bx$)" denote the number of nonzeros in $x$ and $Bx$ obtained by S<small>SNAL</small> using the following estimation

$$\mathrm{nnz}(y) := \min\left\{k \mid \sum_{i=1}^{k} |\hat{y}_i| \geq 0.999\|y\|_1\right\},$$

where $\hat{y}$ is obtained by sorting $y$ such that $|\hat{y}_1| \geq \ldots \geq |\hat{y}_n|$. As can be observed, S<small>SNAL</small> is much faster than all the other four first-order methods. For example, for the `ovarianP` instances, S<small>SNAL</small> can be over 60 times faster than ADMM and 220 times faster than iADMM. For most of the test instances corresponding to the largest problem `ovarianS`, S<small>SNAL</small> only needs about 20 seconds to solve the problems while the other four algorithms run for 20,000 iterations and take about 15 to 90 minutes to only produce rather inaccurate solutions. Here, the poor performance of SLEP, ADMM, iADMM and LADMM indicates that these first-order methods are incapable of obtaining reasonably accurate solutions for difficult large-scale problems.

Comparing ADMM and iADMM, we note that since the sample sizes $m$ for all the tested problems are relatively small (less than 300), ADMM is generally faster as the average cost of solving the $m \times m$ linear system in each iteration is cheaper for ADMM. But we shall see in the next subsection that iADMM would be faster than ADMM when solving problems with large $m$.

Figure 1 presents the performance profiles of S<small>SNAL</small>, ADMM, iADMM and LADMM for all the 80 tested problems. SLEP is not included since it fails on all the test instances. Recall that a point $(x, y)$ is in the performance profile curve of a method if and only if it can solve exactly $(100y\%)$ of all the tested instances in at most $x$ times of the best method for each instance. It can be seen that SSNAL outperforms the other 3 methods by a very large margin. Indeed, SSNAL is more than 20 times faster than all the other tested algorithms for the vast majority of the tested instances.

Table 1: The performance of various algorithms on fused lasso regularized least squares problems on high-dimensional biomedical datasets (accuracy $\eta \leq 10^{-6}$). $m$ is the sample size and $n$ is the dimension of features. In the table, "a" = SSNAL, "b" = SLEP, "c" = ADMM, "d" = iADMM, and "e" = LADMM. "nnz" denotes the number of nonzeros in the solution obtained by SSNAL.

| probname | $\alpha_1; \alpha_2$ | nnz($x$) ; nnz($Bx$) | $\eta$ a \| b \| c \| d \| e | time (hours:minutes:seconds) a \| b \| c \| d \| e |
|---|---|---|---|---|
| DLBCLN | $10^{-3}$ ; 2 | 818 ; 261 | 3.6-7 \| 1.2-5 \| 9.9-7 \| 9.9-7 \| 9.9-7 | 01 \| 16 \| 10 \| 22 \| 18 |
| 160;7399 | $10^{-3}$ ; 0.01 | 157 ; 306 | 9.1-8 \| 4.7-5 \| 9.8-7 \| 6.5-7 \| 3.7-6 | 00 \| 14 \| 17 \| 37 \| 18 |
| $\|A\|_2 = 28.9$ | $10^{-4}$ ; 2 | 848 ; 275 | 3.6-7 \| 3.4-5 \| 9.1-7 \| 8.7-7 \| 3.4-6 | 01 \| 16 \| 20 \| 43 \| 20 |
| | $10^{-4}$ ; 0.01 | 158 ; 306 | 1.5-7 \| 1.0-4 \| 4.3-6 \| 9.9-7 \| 5.8-5 | 01 \| 14 \| 30 \| 1:25 \| 18 |
| lungH1 | $10^{-3}$ ; 2 | 514 ; 325 | 4.9-7 \| 1.1-4 \| 8.7-7 \| 9.3-7 \| 3.9-5 | 01 \| 28 \| 13 \| 22 \| 33 |
| 203;12600 | $10^{-3}$ ; 0.01 | 188 ; 365 | 2.1-7 \| 5.7-4 \| 9.9-7 \| 8.3-7 \| 2.9-2 | 01 \| 24 \| 16 \| 38 \| 29 |
| $\|A\|_2 = 81.5$ | $10^{-4}$ ; 2 | 551 ; 344 | 9.2-8 \| 7.2-4 \| 9.9-7 \| 6.5-7 \| 6.7-2 | 01 \| 28 \| 28 \| 1:02 \| 34 |
| | $10^{-4}$ ; 0.01 | 195 ; 375 | 5.6-8 \| 1.7-3 \| 9.9-7 \| 8.5-7 \| 4.2-2 | 01 \| 25 \| 37 \| 1:34 \| 31 |
| lungH2 | $10^{-3}$ ; 2 | 646 ; 186 | 6.6-8 \| 3.9-5 \| 9.9-7 \| 8.5-7 \| 9.9-7 | 00 \| 26 \| 06 \| 10 \| 19 |
| 149;12533 | $10^{-3}$ ; 0.01 | 137 ; 268 | 4.6-7 \| 1.5-4 \| 9.9-7 \| 7.9-7 \| 2.2-5 | 00 \| 22 \| 24 \| 50 \| 29 |
| $\|A\|_2 = 83.4$ | $10^{-4}$ ; 2 | 775 ; 236 | 2.6-7 \| 1.3-4 \| 8.2-7 \| 9.9-7 \| 2.4-2 | 01 \| 27 \| 15 \| 36 \| 34 |
| | $10^{-4}$ ; 0.01 | 146 ; 285 | 1.2-7 \| 9.9-4 \| 1.4-7 \| 8.3-7 \| 2.0-2 | 01 \| 23 \| 50 \| 1:46 \| 30 |
| ovarianP | $10^{-3}$ ; 2 | 824 ; 144 | 1.6-7 \| 1.6-4 \| 9.9-7 \| 9.0-7 \| 9.7-4 | 01 \| 1:01 \| 19 \| 36 \| 1:06 |
| 253;15153 | $10^{-3}$ ; 0.01 | 180 ; 285 | 1.3-7 \| 6.2-4 \| 9.9-7 \| 9.9-7 \| 2.7-3 | 01 \| 53 \| 59 \| 2:55 \| 1:00 |
| $\|A\|_2 = 114$ | $10^{-4}$ ; 2 | 1259 ; 350 | 2.7-7 \| 4.5-4 \| 8.2-7 \| 9.7-7 \| 1.1-2 | 01 \| 58 \| 45 \| 1:51 \| 1:05 |
| | $10^{-4}$ ; 0.01 | 255 ; 412 | 9.9-8 \| 1.4-3 \| 2.9-5 \| 1.9-5 \| 4.4-3 | 01 \| 55 \| 1:55 \| 6:11 \| 1:01 |
| ovarianS | $10^{-3}$ ; 2 | 1958 ; 352 | 6.7-7 \| 3.8-3 \| 2.1-6 \| 8.5-7 \| 9.4-3 | 15 \| 20:17 \| 46:31 \| 57:19 \| 23:23 |
| 216;373401 | $10^{-3}$ ; 0.01 | 205 ; 409 | 6.3-7 \| 8.5-3 \| 1.6-3 \| 3.6-4 \| 3.8-2 | 14 \| 18:41 \| 41:12 \| 1:15:03 \| 21:55 |
| $\|A\|_2 = 539$ | $10^{-4}$ ; 2 | 1963 ; 380 | 2.5-7 \| 7.3-3 \| 1.1-3 \| 1.8-3 \| 9.0-2 | 20 \| 16:39 \| 45:03 \| 1:17:03 \| 22:59 |
| | $10^{-4}$ ; 0.01 | 212 ; 422 | 2.5-7 \| 6.6-3 \| 1.2-3 \| 6.8-2 \| 1.8-1 | 18 \| 18:15 \| 44:24 \| 1:50:34 \| 23:13 |

### 5.1.2 Numerical results for UCI datasets

In this subsection, we test all the algorithms on the same large scale UCI datasets $(A, b)$ as in [21] that are originally obtained from the LIBSVM datasets [6].

In Table 2, we report the detailed numerical results for SSNAL, SLEP, ADMM, iADMM and LADMM in solving the least squares fused lasso regularized problem (2) on large-scale UCI datasets. In these tests, we choose the regularized parameter $\alpha_2 \in \{1, 0.5, 0.2, 0.01\}$. Meanwhile, in order to produce reasonable non-zeros in the optimal solution $x$ and $Bx$, in these tests, we choose $\alpha_1 \in \{10^{-6}, 10^{-7}\}$ for problems `E2006.train` and `E2006.test`, $\alpha_1 \in \{10^{-5}, 10^{-6}\}$ for problem `bodyfat7` and $\alpha_1 \in \{10^{-3}, 10^{-4}\}$ for all the other instances. In total, we tested 80 instances. Note that in order to save space, Table 2 only reports the results for a subset of these instances. We also present in Figure 2 the performance profiles of SSNAL, SELP, ADMM, iADMM and LADMM for all the tested problems.

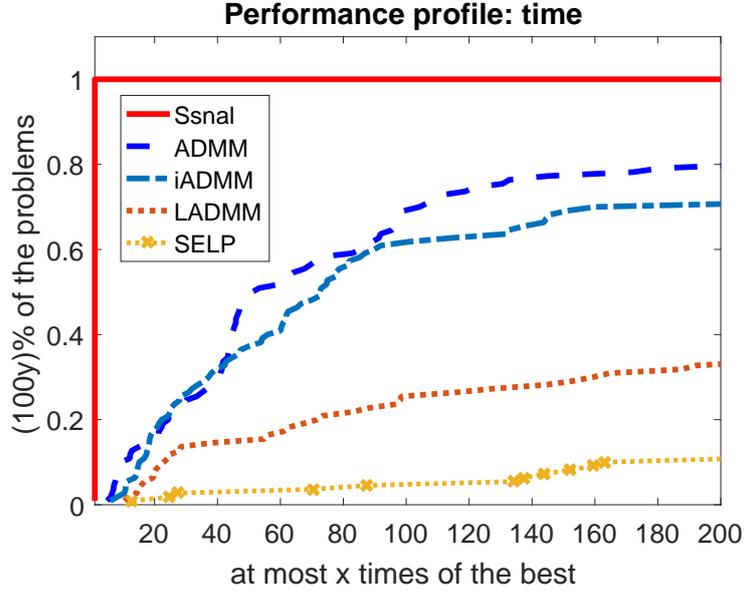

Figure 2: Performance profiles of SSNCG, SELP, ADMM, iADMM and LADMM on UCI datasets.

Table 2: Same as Table 1 but for large-scale UCI datasets.

| probname $m;n$ | $\alpha_1;\alpha_2$ | nnz(x) ; nnz(Bx) | $\eta$ a \| b \| c \| d \| e | time (hours:minutes:seconds) a \| b \| c \| d \| e |
|---|---|---|---|---|
| E2006.train 16087;150360 $\|A\|_2 = 437$ | $10^{-6}$ ; 0.5 | 8 ; 13 | 2.1-7 \| 2.1-3 \| 2.5-7 \| 3.4-7 \| 3.2-3 | 03 \| 18:49 \| 36:42 \| 7:34 \| 18:39 |
| | $10^{-6}$ ; 0.01 | 25 ; 47 | 6.1-8 \| 9.9-4 \| 8.8-8 \| 5.3-7 \| 6.4-3 | 04 \| 21:12 \| 50:04 \| 9:48 \| 19:48 |
| | $10^{-7}$ ; 0.5 | 657 ; 1069 | 9.3-7 \| 4.2-3 \| 2.7-7 \| 4.9-8 \| 8.5-3 | 19 \| 20:13 \| 42:04 \| 13:29 \| 20:05 |
| | $10^{-7}$ ; 0.01 | 1424 ; 2764 | 1.7-7 \| 4.3-3 \| 4.2-7 \| 8.9-7 \| 9.1-3 | 1:13 \| 19:41 \| 45:20 \| 18:01 \| 20:14 |
| E2006.test 3308;150358 $\|A\|_2 = 219$ | $10^{-6}$ ; 0.5 | 14 ; 24 | 2.6-8 \| 6.8-4 \| 3.6-8 \| 5.7-7 \| 4.2-3 | 02 \| 5:32 \| 2:55 \| 2:17 \| 6:37 |
| | $10^{-6}$ ; 0.01 | 49 ; 95 | 1.7-8 \| 4.9-4 \| 9.2-8 \| 3.3-7 \| 4.7-3 | 02 \| 5:20 \| 2:59 \| 2:22 \| 6:39 |
| | $10^{-7}$ ; 0.5 | 765 ; 1384 | 2.8-8 \| 1.2-3 \| 5.4-7 \| 3.9-7 \| 5.1-3 | 12 \| 5:31 \| 4:10 \| 4:56 \| 6:53 |
| | $10^{-7}$ ; 0.01 | 1317 ; 2581 | 2.9-7 \| 1.1-3 \| 7.1-7 \| 7.8-7 \| 5.1-3 | 53 \| 5:07 \| 4:23 \| 5:15 \| 6:24 |
| log1p.E2006.train 16087;4272227 $\|A\|_2 = 7650$ | $10^{-3}$ ; 0.5 | 4 ; 5 | 4.0-8 \| 2.5-4 \| 3.5-7 \| 2.7-7 \| 5.3-3 | 24 \| 2:52:08 \| 36:17 \| 14:33 \| 3:00:01 |
| | $10^{-3}$ ; 0.01 | 5 ; 9 | 9.6-8 \| 1.8-5 \| 6.5-7 \| 4.7-7 \| 2.4-3 | 25 \| 2:45:56 \| 43:08 \| 16:21 \| 3:00:01 |
| | $10^{-4}$ ; 0.5 | 256 ; 340 | 1.2-7 \| 1.3-4 \| 9.9-7 \| 8.8-7 \| 1.3-2 | 53 \| 2:47:20 \| 52:33 \| 32:45 \| 3:00:01 |
| | $10^{-4}$ ; 0.01 | 576 ; 1100 | 9.8-7 \| 1.6-4 \| 7.8-7 \| 6.7-7 \| 1.4-2 | 1:09 \| 2:44:02 \| 1:01:16 \| 54:23 \| 3:00:01 |
| log1p.E2006.test 3308;4272226 $\|A\|_2 = 3830$ | $10^{-3}$ ; 0.5 | 4 ; 5 | 6.1-7 \| 1.5-4 \| 1.5-8 \| 5.3-7 \| 8.3-4 | 17 \| 1:44:15 \| 6:20 \| 6:05 \| 1:56:52 |
| | $10^{-3}$ ; 0.01 | 8 ; 15 | 5.8-8 \| 1.4-4 \| 2.7-7 \| 2.0-7 \| 5.2-4 | 21 \| 1:39:06 \| 8:32 \| 8:45 \| 1:52:15 |
| | $10^{-4}$ ; 0.5 | 597 ; 842 | 7.3-8 \| 2.5-4 \| 5.1-7 \| 6.8-7 \| 1.6-3 | 58 \| 1:41:13 \| 11:54 \| 14:40 \| 1:54:20 |
| | $10^{-4}$ ; 0.01 | 1059 ; 2035 | 2.0-7 \| 2.2-4 \| 2.5-7 \| 9.8-7 \| 2.7-3 | 42 \| 1:35:41 \| 12:16 \| 13:37 \| 1:46:12 |
| pyrim5 74;201376 $\|A\|_2 = 1110$ | $10^{-3}$ ; 0.5 | 174 ; 123 | 8.5-7 \| 5.6-3 \| 9.9-7 \| 9.9-7 \| 3.2-4 | 04 \| 8:27 \| 12:17 \| 29:20 \| 9:58 |
| | $10^{-3}$ ; 0.01 | 75 ; 145 | 1.7-7 \| 1.9-3 \| 6.8-5 \| 2.0-4 \| 4.3-4 | 04 \| 7:33 \| 16:17 \| 59:36 \| 8:31 |
| | $10^{-4}$ ; 0.5 | 233 ; 142 | 3.0-7 \| 6.8-3 \| 6.2-5 \| 1.4-4 \| 2.5-3 | 07 \| 8:27 \| 19:55 \| 1:35:54 \| 9:33 |
| | $10^{-4}$ ; 0.01 | 91 ; 156 | 3.1-8 \| 6.2-3 \| 7.1-3 \| 1.7-3 \| 1.6-3 | 06 \| 7:26 \| 17:16 \| 1:42:19 \| 8:21 |
| triazines4 186;635376 $\|A\|_2 = 4550$ | $10^{-3}$ ; 0.5 | 679 ; 260 | 2.9-7 \| 3.4-3 \| 2.9-5 \| 5.9-3 \| 3.8-3 | 25 \| 1:02:49 \| 2:07:21 \| 3:00:02 \| 1:03:59 |
| | $10^{-3}$ ; 0.01 | 217 ; 302 | 1.4-7 \| 2.8-3 \| 3.0-4 \| 1.4-1 \| 4.8-3 | 27 \| 54:52 \| 1:56:48 \| 3:00:01 \| 56:11 |
| | $10^{-4}$ ; 0.5 | 875 ; 334 | 4.5-7 \| 1.2-2 \| 9.9-3 \| 8.6-1 \| 5.2-2 | 37 \| 1:00:20 \| 2:13:51 \| 3:00:10 \| 1:02:02 |
| | $10^{-4}$ ; 0.01 | 223 ; 355 | 3.3-7 \| 1.3-2 \| 2.6-2 \| 7.8-1 \| 2.3-2 | 40 \| 1:07:27 \| 2:37:41 \| 3:00:01 \| 1:11:24 |
| bodyfat7 | $10^{-5}$ ; 0.5 | 36 ; 29 | 6.4-7 \| 2.1-4 \| 9.9-7 \| 6.0-7 \| 2.6-6 | 03 \| 7:14 \| 2:22 \| 7:59 \| 7:39 |

Table 2: Same as Table 1 but for large-scale UCI datasets.

| probname $m;n$ | $\alpha_1;\alpha_2$ | nnz($x$) ; nnz($Bx$) | $\eta$ a \| b \| c \| d \| e | time (hours:minutes:seconds) a \| b \| c \| d \| e |
|---|---|---|---|---|
| 252;116280 | $10^{-5}$ ; 0.01 | 25 ; 43 | 3.1-7 \| 7.8-4 \| 5.8-7 \| 8.4-7 \| 2.8-5 | 03 \| 6:56 \| 2:19 \| 8:59 \| 7:20 |
| $\|A\|_2 = 230$ | $10^{-6}$ ; 0.5 | 142 ; 136 | 7.9-7 \| 1.4-3 \| 9.2-7 \| 9.4-7 \| 9.1-4 | 05 \| 7:05 \| 3:03 \| 19:30 \| 7:33 |
| | $10^{-6}$ ; 0.01 | 101 ; 190 | 9.2-8 \| 1.4-3 \| 9.9-7 \| 9.9-7 \| 8.7-4 | 06 \| 6:41 \| 9:33 \| 47:33 \| 7:02 |
| housing7 | $10^{-3}$ ; 0.5 | 126 ; 149 | 4.0-7 \| 2.1-4 \| 9.9-7 \| 9.2-7 \| 2.4-6 | 02 \| 7:21 \| 4:20 \| 17:24 \| 7:33 |
| 506;77520 | $10^{-3}$ ; 0.01 | 151 ; 284 | 4.8-7 \| 3.6-4 \| 9.9-7 \| 7.6-7 \| 1.0-4 | 02 \| 7:03 \| 4:26 \| 18:11 \| 7:16 |
| $\|A\|_2 = 573$ | $10^{-4}$ ; 0.5 | 253 ; 352 | 1.6-7 \| 4.4-3 \| 9.9-7 \| 7.6-7 \| 4.0-4 | 03 \| 7:26 \| 6:44 \| 1:02:18 \| 7:36 |
| | $10^{-4}$ ; 0.01 | 276 ; 543 | 2.7-7 \| 7.3-3 \| 9.9-7 \| 7.0-7 \| 1.3-3 | 04 \| 6:59 \| 8:55 \| 1:36:21 \| 7:16 |

It can be clearly observed in Table 2 and Figure 2 that Ssnal outperforms all the other tested first-order algorithms by a large margin where the factor can be up to at least 150 times faster. In fact, for over 80 percent of the instances, Ssnal is at least 20 times faster than all the other tested algorithms. We also note that Ssnal is the only algorithm which can solve all the test instances to the required accuracy. For the test instances corresponding to problem `triazines4`, Ssnal only needs less than 1 minute to produce a solution with the required accuracy while all the other first-order algorithms spend over 1 hour (2 and 3 hours for ADMM and iADMM) to only produce poor accuracy solutions with $\eta \approx 10^{-3}$. These observations again demonstrate the power of Ssnal over the other tested first-order algorithms. Moreover, from the unfavorable performance of SLEP, ADMM, iADMM and LADMM, one can safely conclude that these first-order algorithms can only be used to solve relatively easy problems. This fact together with the superior efficiency and robustness of Ssnal indicates that it is necessary to incorporate second-order nonsmooth analysis into the algorithmic design, especially when solving large scale difficult problems. In particular, the efficiency of our Ssnal depends critically on our ability to extract and exploit the underlying second-order structured sparsity in the problems.

Among the first-order methods, one can observe from the results corresponding to `E2006.train` and `log1p.E2006.train` that when $m$ (sample size) is large, iADMM performs better than ADMM. This demonstrates the advantage of iADMM over ADMM in having the flexibility of solving large $m \times m$ linear systems in the subproblems inexactly by an iterative solver such as the PCG method. On the other hand, when $m$ is relatively small, it is advantageous to solve the $m \times m$ linear systems in the subproblems of ADMM by a direct solver as opposed to using an iterative solver as in the case of iADMM.

## 5.2 Numerical results for least squares constrained fused lasso problems

Given $\delta > 0$, recall the least squares constrained fused lasso problems given in (1)

$$\min \{p(x) \mid \|Ax - b\| \leq \delta\}. \tag{32}$$

In this section, we present the numerical results obtained by a level-set method for solving the problem. Since Ssnal is applied to solve the regularized least squares subproblems (33), we term the algorithm as the Ssnal based level-set method (in short, Ssnal-LSM). More specifically, Ssnal-LSM is based on a bisection method to solve the univariate nonlinear equation associated with the value function $\varphi$ given in (4)

$$\varphi(\mu) = \delta.$$

At the $k$-th iteration of SSNAL-LSM, $\varphi(\mu_k)$ is evaluated through using SSNAL to solve the subproblem (33) for the given parameter $\mu_k \geq 0$.

---

**Algorithm SSNAL-LSM: An SSNAL based level-set method for (32).**

Let $\mu_\infty > \mu_0 \geq 0$ be two given parameters. Set $\underline{\mu} = \mu_0$, $\overline{\mu} = \mu_\infty$ and $\mu_1 = (\underline{\mu} + \overline{\mu})/2$. For $k = 1, 2, \ldots$, perform the following steps in each iteration:

**Step 1.** Use SSNAL to compute
$$x^k = \arg\min \left\{ \frac{1}{2}\|Ax - b\|^2 + \mu_k p(x) \right\}. \tag{33}$$

**Step 2.** Compute $\varphi(\mu_k) = \|Ax^k - b\|$. If $\varphi(\mu_k) > \delta$, update $\overline{\mu} = \mu_k$, otherwise, $\underline{\mu} = \mu_k$.

**Step 3.** Update $\mu_{k+1} = (\underline{\mu} + \overline{\mu})/2$.

---

For testing purpose, the fused lasso regularizer $p$ is chosen as follows
$$p(x) = \|x\|_1 + 2\|Bx\|_1, \quad \forall x \in \Re^n.$$

The noise level controlling parameter $\delta$ in (32) is chosen to be $\delta = \gamma\|b\|$, where $0 < \gamma < 1$. We choose the initial parameters $\mu_0 = 0$ and $\mu_\infty = \|A^T b\|_\infty$. In our numerical experiments, we measure the accuracy of an approximate optimal solution $\tilde{x}$ for (32) by using the following relative residual:
$$\eta = \frac{|\tilde{\varphi} - \delta|}{\max\{1, \delta\}},$$
where $\tilde{\varphi} := \|A\tilde{x} - b\|$. We stop the algorithm when $\eta \leq 10^{-6}$. In solving the subproblems (33) by the SSNAL method, the required accuracy for $x^k$ is set to $10^{-8}$. The large scale UCI and biomedical datasets are both used in the experiments.

In Table 3, we report the detailed results for SSNAL-LSM in solving the least squares constrained fused lasso problems of form (32) for large-scale UCI and biomedical datasets. In our tests, we choose three different $\gamma$ for each test instance to show the changes in the sparsity patterns of the obtained solutions. In the table, $\mu^*$ denotes the solution for $\varphi(\mu) = \delta$ for a given $\delta$. The column "iteration" reports the number of iterations taken by the SSNAL-LSM to solve the problems. It can be seen from the table that SSNAL-LSM usually takes about 20 to 30 iterations to achieve a sparse solution with the desired accuracy. That is, we only need to use SSNAL to solve 20 to 30 regularized least squares fused lasso subproblems. Combining the superior performance of SSNAL presented in Subsection 5.1, one can safely conclude that for most test instances, the time required by SSNAL-LSM to solve the constrained problem (32) can still be much less than that required by any of the previously tested first-order methods to solve a single fused lasso regularized least squares problem.

## 6 Conclusion

In this paper, we showed that the level-set method can be used to solve least squares constrained fused lasso problems where the subproblems are fused lasso regularized least squares problems. As the backbone of the level-set method, we designed an extremely fast semismooth Newton based augmented Lagrangian method, i.e., SSNAL, for solving the fused lasso regularized least squares problems. We achieve the superior performance of SSNAL through a careful analysis of the structures of the generalized Jacobian for the proximal mapping of the fused lasso regularizer. In particular, we uncovered crucial second-order structured sparsity in the used generalized Jacobian and designed several delicate numerical techniques to exploit the underlying structures for solving the semismooth Newton systems in the SSNAL algorithm very efficiently. Extensive numerical experiments on fused lasso regularized least squares problems on high-dimensional real data instances show the great benefits of our second-order nonsmooth analysis based algorithms.

Table 3: The performance of SSNAL-LSM on least squares constrained fused lasso problem (32) on large-scale datasets (accuracy $\eta \leq 10^{-6}$). $m$ is the sample size and $n$ is the dimension of features. "nnz" denotes the number of nonzeros in the solution. The computation time is in the format of "hours:minutes:seconds".

| probname | $\gamma$ | nnz($x$) ; nnz($Bx$) | $\mu^*$ | iteration | $\eta$ | time |
| $m;n$ | | | | | | |
|---|---|---|---|---|---|---|
| E2006.train | 1.0-1 | 840 ; 547 | $1.30 - 2$ | 40 | 1.2-7 | 3:07 |
| 16087;150360 | 1.5-1 | 1 ; 1 | $6.86 + 3$ | 22 | 1.5-7 | 13 |
| | 2.0-1 | 1 ; 1 | $1.11 + 4$ | 22 | 2.0-7 | 12 |
| log1p.E2006.train | 1.0-1 | 345 ; 177 | $2.38 + 1$ | 27 | 1.3-7 | 9:17 |
| 16087;4272227 | 1.5-1 | 20 ; 6 | $2.49 + 3$ | 25 | 7.0-7 | 4:55 |
| | 2.0-1 | 20 ; 6 | $3.94 + 3$ | 25 | 3.3-7 | 4:46 |
| E2006.test | 5.0-2 | 2393 ; 2240 | $1.20 - 3$ | 43 | 5.9-7 | 4:45 |
| 3308;150358 | 7.5-2 | 603 ; 680 | $3.36 - 3$ | 41 | 5.1-8 | 1:08 |
| | 1.0-1 | 1 ; 1 | $2.56 + 2$ | 21 | 5.8-7 | 06 |
| log1p.E2006.test | 5.0-2 | 3685 ; 2609 | $1.40 + 0$ | 34 | 2.3-7 | 15:40 |
| 3308;4272226 | 7.5-2 | 1504 ; 1003 | $3.27 + 0$ | 31 | 5.8-7 | 8:59 |
| | 1.0-1 | 20 ; 7 | $1.15 + 2$ | 24 | 8.5-7 | 4:10 |
| pyrim5 | 1.0-1 | 254 ; 49 | $3.74 - 1$ | 24 | 4.7-8 | 40 |
| 74;201376 | 2.0-1 | 38 ; 9 | $1.18 + 0$ | 24 | 5.6-7 | 37 |
| | 3.0-1 | 54 ; 10 | $2.45 + 0$ | 23 | 7.1-8 | 35 |
| triazines4 | 1.0-1 | 1338 ; 194 | $1.59 - 1$ | 27 | 8.2-7 | 5:32 |
| 186;635376 | 2.0-1 | 782 ; 47 | $1.91 + 0$ | 23 | 9.4-7 | 3:26 |
| | 3.0-1 | 243 ; 18 | $9.33 + 0$ | 19 | 4.7-7 | 2:35 |
| housing7 | 1.0-1 | 238 ; 134 | $3.87 + 0$ | 28 | 5.7-7 | 36 |
| 506;77520 | 2.0-1 | 34 ; 21 | $7.06 + 1$ | 23 | 5.0-7 | 20 |
| | 3.0-1 | 17 ; 12 | $1.83 + 2$ | 21 | 4.0-7 | 17 |
| bodyfat7 | 1.0-4 | 731 ; 391 | $2.70 - 6$ | 35 | 4.9-8 | 1:00 |
| 252;116280 | 1.0-3 | 322 ; 150 | $9.87 - 5$ | 33 | 7.0-7 | 43 |
| | 1.0-2 | 2 ; 3 | $1.98 - 1$ | 25 | 4.1-7 | 19 |
| ovarianP | 1.5-1 | 591 ; 53 | $5.66 - 2$ | 25 | 1.1-7 | 06 |
| 253;15153 | 2.0-1 | 686 ; 30 | $1.53 - 1$ | 22 | 1.1-9 | 05 |
| | 2.5-1 | 368 ; 22 | $2.65 - 1$ | 23 | 7.0-7 | 05 |
| ovarianS | 1.5-1 | 1506 ; 218 | $6.02 - 2$ | 26 | 4.6-7 | 2:05 |
| 216;373401 | 2.0-1 | 1395 ; 175 | $8.40 - 2$ | 25 | 4.2-7 | 1:53 |
| | 2.5-1 | 1123 ; 133 | $1.11 - 1$ | 23 | 5.5-7 | 1:43 |